\documentclass[10pt]{article}

\usepackage{amsfonts}
\usepackage{amsmath}
\usepackage{graphicx}
\usepackage{psfrag}
\usepackage{CJK}

\usepackage[centerlast]{caption}
 \usepackage{subfig}

\newtheorem{satz}{Satz}[section]
\newtheorem{theorem}[satz]{\indent\textsc Theorem}
\newtheorem{lemma}[satz]{Lemma}

\newtheorem{remark}[satz]{Remark}
\newtheorem{definition}[satz]{Definition}

\def\proofbegin{{\textit  Proof}. }
\allowdisplaybreaks

\title{Preservation of Takens-Bogdanov bifurcations for delay differential equations by Euler discretization\thanks{Supported by NSFC grants 10971022 and 11071102.}} 

\author{{\sc Yingxiang Xu}\footnote{School of Mathematics and Statistics, Northeast Normal University, Changchun 130024, China.
E-mail: yxxu@nenu.edu.cn} and {\sc Chengchun Gong}\footnote{School of Mathematics, Jilin University, Changchun 130012, China. E-mail: gongcc@jlu.edu.cn}}
\date{}%
\begin{document}%
\maketitle

\medskip

\medskip
%

%

\begin{abstract}
A new technique for calculating the normal forms associated with the map restricted to the center manifold of a class of parameterized maps near the fixed point is given first. Then we show the Takens-Bogdanov point of delay differential equations is inherited by the forward Euler method without any shift and turns into a 1:1 resonance point. The normal form near the 1:1 resonance point for the numerical discretization is calculated next by applying the new technique to the map defined by the forward Euler method. The local dynamical behaviors are analyzed in detail through the normal form. It shows the Hopf point branch and the homoclinic branch emanating from the Takens-Bogdanov point are $O(\varepsilon)$ shifted by the forward Euler method, where $\varepsilon$ is step size. At last, a numerical experiment is carried to show the results.
\end{abstract}

\section{Introduction}

Numerical methods may take on many possible phenomena when applied to certain dynamical systems. It is of great importance to investigate what kind of properties of the original systems could be preserved by discretization. Tremendous researches on numerical stabilities of sorts are the first investigation of this type where various sufficient and necessary conditions are developed for numerical schemes  reproducing the asymptotic stability of differential equations, one can refer to the monographs \cite{Bu08} and \cite{Bz03} for ordinary  differential equations (ODEs) and delay differential equations (DDEs), respectively. Whether or not the numerical discretization will inherit the bifurcations of the original systems is another important research field of this type.

In this paper, we consider the DDEs of the type
\begin{equation}\label{intro1}
\dot{z}(t)=f(z(t),z(t-1),\alpha),
\end{equation}
where $z\in \mathbb{R}^n$, $\alpha\in\mathbb{R}^2$ is a bifurcation parameter, $f(z_1,z_2,\alpha)$ is $C^r(r\geq 2)$ smooth with respect to $z_1,z_2$ and $\alpha$. The state space of (\ref{intro1}), denoted by $C=C([-1,0],\mathbb{R}^n)$, is a Banach space of continuous mappings from $[-1,0]$ to $\mathbb{R}^n$ equipped with norm $\|\phi\|=\max_{\theta\in[-1,0]}|\phi(\theta)|$ ($|\cdot|$ is some norm in $\mathbb{R}^n$). The equation (\ref{intro1}) is assumed to undergo a Takens-Bogdanov bifurcation near $(z,\alpha)=(0,0)$. Then, in the parameter plane $(\alpha_1,\alpha_2)$ there exist a Hopf point branch and a homoclinic branch emanating from the Takens-Bogdanov point of equation (\ref{intro1}), see \cite{xu08}. In this paper we show that the forward Euler discretization, when applied to equation (\ref{intro1}), can preserve the bifurcation structure near the Takens-Bogdanov point $(z,\alpha)=(0,0)$ of (\ref{intro1}) by an $O(\varepsilon)$ shift, where $\varepsilon$ is the step size of the Euler method.

Problems of this type have been focused on ODEs for many years. The Hopf bifurcation accepted the most attentions in this field after the work of Hofbauer and Iooss \cite{Hi84}, where they proved the forward Euler discretization exhibits the Hopf bifurcation of the same type as the continuous system undergoes. In 1998, Wang, Blum and Li \cite{Wbl98} gave the most complete results for the codimension 1 bifurcations: the general one step methods of order $p$ and specific methods like Euler, backward Euler, explicit and implicit Runge-Kutta methods are proved to inherit the elementary bifurcations, such as saddle-node, fold, pitchfork, Hopf and transcritical bifurcations, of the continuous systems. For connecting orbits, Beyn \cite{be87} first showed the existence of homoclinic orbits under discretization. Fiedler and Scheurle \cite{fs96}, Zou and Beyn \cite{Zb03} proved the general one step method exhibits the discrete connecting orbits approximating to the one of the original ODEs by the order of the method, independently.  
Most recently, the major concern for this problem is whether the numerical discretization can inherit the codimension 2 bifurcations of ODEs. L\'oczi and Ch\'avez \cite{Lc09} showed the Runge-Kutta method will reproduce the nondegenerate fold, cusp and Takens-Bogdanov singularities of general ODEs without any shift. For generalized Hopf bifurcations, the Hopf point branch emanating from the generalized Hopf point is shifted by the order of the general one step method used as expected, however, the generalized Hopf point is shifted of the first order, regardless the order of the numerical method used, and turns into generalized Neimark-Sacker points. This result seems not natural, but the numerical experiments shows a much better result cannot be  expected, cf. {\cite{Ch11}}. For the Hopf-Hopf bifurcation, the one step method will reproduce the bifurcation point as a double Neimark-Sacker point by an $O(\varepsilon^p)$ shift, as well as the Hopf point branch bifurcated from the point, see \cite{Ch12}. Ch\'avez \cite{Ch10}, under the assumption that the ODEs undergo the fold-Hopf or Takens-Bogdanov bifurcation, proved the singularities will be persisted by the Runge-Kutta method and the Hopf point branch emanating from the singular points is $O(\varepsilon^p)$ shifted ($\varepsilon$ is step size). In addition,  the discrete fold point branch and the discrete Hopf point branch emanating from the singularities is shown to intersect transversally. However, according to what we have seen up to now, no discussion about the preservation of the homoclinic branch emanating from the Takens-Bogdanov point by the numerical scheme is available.

In these days, the preservation of dynamical behaviors for DDEs by numerical discretizations received many considerations as well. Wulf and Ford first showed the forward Euler method will reproduce the same type Hopf bifurcations as the scalar DDE with one constant delay undergoes \cite{Wf00}. The subsequent numerous relative works had led this problem to a more general case of both the numerical methods (Runge-Kutta method, linear multi-step method etc.) and the type of the differential equations (with multiple delays or state-dependent delay, etc.). For other dynamics of DDEs preserved by discretization, we recall the following results. Liu, Gao and Yang \cite{Lgy09} proved the Runge-Kutta method preserves the oscillations of the equation $\dot x(t)+ax(t)+a_1x([t-1])=0$. The asymptotically stable periodic orbits of the autonomous DDEs with one constant delay were proved by In't Hout and Lubich \cite{hl98} to be shifted by the $s$-stage Runge-Kutta methods  with the discretization order $p$. For a more general delay equation
\begin{equation}\label{intro5}
\dot x=f(x_t),
\end{equation}
where $x_t(\theta)=x(t+\theta)$, Farkas proved by explicit calculation that the unstable manifold will be close to its Euler discretized counterpart if the step size $\varepsilon$ is sufficiently small \cite{f01} as well as a numerical shadowing result which obviously means the stable manifolds of (\ref{intro5}) could be preserved numerically \cite{f02}. Xu and Zou \cite{Xz11} extended the results of \cite{Zb03} to DDEs, they showed the homoclinic orbits should be preserved by the forward Euler method with a shift of $O(\varepsilon)$. In addition, the possibility of extending the numerical scheme to the implicit method and the Runge-Kutta method is discussed there.

The problem we care most is whether the bifurcation diagram near Takens-Bogdanov point of DDEs could be inherited by numerical discretiztation. This requires us to investigate the dynamical behaviors of the map defined by the numerical scheme. The normal form analysis is a useful tool for accomplishing this aim. It allows us to know how the normal form coefficients and the generalized eigenspace of the numerical scheme are related to their continuous counterparts. Consequently, the dynamics of the numerical scheme could be studied in detail in terms of the coefficients of the original equation and then we know how the bifurcation diagram is preserved by direct comparison. Trying to fulfill our goal in this way, we will also face some problems.

Generally speaking, the numerical scheme for ODEs could be regarded as a map. For example, the one step method for solving
$$
\dot x=f(x),\ x\in \mathbb{R}^n
$$
has the general form
\begin{equation}\label{numsch}
x_{k+1}=x_{k}+\varepsilon\phi(x_{k},\varepsilon),
\end{equation}
where $\varepsilon$ is step size. Obviously, equation (\ref{numsch}) is a map from $\mathbb{R}^n$ to $\mathbb{R}^n$ and the method of normal forms for finite-dimensional ODEs could be modified to analyze the local bifurcation behavior. Where the representations of the complete eigenvectors are required to construct a transformation which lead the map (\ref{numsch}) to a standard form could be dealt with easily (see \cite{Wi90} section 1.1C for detail). In addition, the technique requires computing the center manifold first before evaluating normal forms for the map on the center manifold. However, the forward Euler method for solving
$$
\dot x(t)=ax(t)+bx(t-1),\ x\in\mathbb{R}^n
$$
has the form
\begin{equation}\label{numsch2}
x_{k+1}=x_k+\varepsilon au_k+\varepsilon bx_{k-m},
\end{equation}
where $\varepsilon=\frac{1}{m}$ is step size, $m$ is a positive integer. The map defined by (\ref{numsch2}) is,  obviously,  from $\mathbb{R}^{(m+1)\times n}$ to $\mathbb{R}^{(m+1)\times n}$ and the dimensionality tends to infinity as the step size tends to 0. As a result, the amount of the eigenvectors is very large as the step size is small, and the method of normal forms for maps described in \cite{Wi90} is not efficient any more.

Highlighted by the normal form analysis for retarded functional differential equations produced by Faria \cite{Faria95a, Fa95b}, a new method for calculating the normal forms associated with the map restricted to the center manifold tangent to the invariant manifold of the lienarization of the map of the type mentioned above at the fixed point is given. It has many advantages: one, without computing beforehand the center manifold of the singularity; two, allows obtaining the coefficients in the normal form explicitly in terms of the considered map, therefore, one can obtain the accurate dynamics through analyzing the obtained normal form, and three, no need of computing the eigenvectors associated with the noncritical eigenvalues, and then the amount of calculations are reduced remarkably.

The rest of the paper is arranged as follows. We first go deep into the techniques for normal form calculations for the maps similar to (\ref{numsch2}) in section 2. The Takens-Bogdanov bifurcations for equation (\ref{intro1}) is recalled in section 3. In Section 4, we first show the Takens-Bogdanov point of DDEs is inherited by the forward Euler method exactly and turns into a 1:1 resonance point. Then the normal form near the 1:1 resonance point for the numerical method is calculated by applying the techniques developed in Section 2. The local dynamical behaviors are analyzed in detail through the obtained normal form. It shows that the Hopf point branch and the homoclinic branch bifurcated from the Takens-Bogdanov point of (\ref{intro1}) are inherited by the forward Euler scheme by a shift of $O(\varepsilon)$. A numerical experiment for a scalar DDE is carried to show our theoretical results in the last section.

\section{Normal forms for a class of maps with parameters}
We consider a class of maps of the following type
\begin{equation}\label{map}
u\mapsto C(\alpha,\varepsilon)u+F(u,\alpha,\varepsilon),
\end{equation}
where $u\in\mathbb{R}^{(m+1)nl}$ with $n,m,l\in \mathbb{Z}_+$, $\varepsilon=\frac{1}{m}$ is a parameter used in later applications where $\varepsilon$ corresponds to the step size, $\alpha=(\alpha_1,\cdots,\alpha_p)\in\mathbb{R}^p$ is a bifurcation parameter. $F$ is $C^r(r\geq2)$ smooth with respect to $u$ and $\alpha$, and satisfies $F(0,0,\varepsilon)=0$ and $\frac{\partial F}{\partial u}(0,0,\varepsilon)=0$. From the map (\ref{map}), we see that the dimensionality of (\ref{map}) will tend to infinity as the parameter $\varepsilon$ tends to $0$.

Denote $C(\varepsilon)=C(0,\varepsilon)$ and $\bar F(u,\alpha,\varepsilon)=C(\alpha,\varepsilon)-C(0,\varepsilon)+F(u,\alpha,\varepsilon)$, then we can write the map (\ref{map}) as
\begin{equation}
u\mapsto C(\varepsilon)u+\bar F(u,\alpha,\varepsilon).
\end{equation}
We drop the bar from $\bar F$ if no confusion occurs. For simplicity, we will drop the dependence on $\varepsilon$ in the rest of this section.

We reformulate the map (\ref{map}) by regarding $\alpha$ as a variable to the following map without parameters
\begin{equation}\label{suspend}
\dbinom{u}{\alpha}\mapsto\tilde C\dbinom{u}{\alpha}+\dbinom{F(u,\alpha)}{0_p},
\end{equation}
where $\tilde C= {\rm diag}\{C,I_p\}$. To linearize the map (\ref{suspend}) at $(u,\alpha)=(0,0)$ we obtain
\begin{equation}\label{linmap}
\dbinom{u}{\alpha}\mapsto\tilde C\dbinom{u}{\alpha},
\end{equation}
which has the characteristic equation given by
\begin{equation}
\det(
\lambda I-\tilde C)=0\Leftrightarrow (\lambda-1)^p\det (\lambda I-C)=0.
\end{equation}
As seen in the former equation, we often omit the subscript of the identity matrix $I$ if no confusion occurs throughout this paper.

%
Denoted by $\Lambda=\{\lambda\in \sigma(C):|\lambda|=1\}$ and $\tilde{\Lambda}=\Lambda\cup\{\underbrace{1,\cdots,1}_{p\ times}\}$,
$\phi_i$ the (generalized) eigenvectors of $C$ associated with eigenvalues $\lambda_i$, $i=1,2,\cdots,(m+1)nl$, then $P={\rm span}\{\phi_i:\lambda_i\in\Lambda\}$ is the invariant space of $C$ associated with $\Lambda$. If $c$ is the number of the eigenvalues of $C$ in $\Lambda$, counting multiplicities, then we have $\dim P=c$.
If we denote by $\Phi_c=(\phi_1,\cdots,\phi_c)$ a basis for $P$, $\Psi_c=(\psi_1^T,\cdots,\psi_c^T)^T$ a basis for the dual space $P^\ast$ in $\mathbb{R}^{(m+1)nl\ast}$, then $(\Psi_c,\Phi_c)=I_c$, the identity matrix in $\mathbb{R}^c$, where the dual form takes the scalar product of vectors. We denote by $J$ the $c\times c$ constant matrix such that $C\Phi_c=\Phi_c J$; its spectrum coincides with $\Lambda$. Denoted by $Q$ the invariant space of $C$ associated with $\sigma(C)\setminus \Lambda$, then $Q={\rm span}\{\phi_i, i=c+1,\cdots,(m+1)nl\}$. If we denote $\Phi_{su}=(\phi_{c+1},\cdots,\phi_{(m+1)nl})$ a basis for $Q$, $\Psi_{su}=(\psi_{c+1}^T,\cdots,\psi_{(m+1)nl}^T)^T$ a basis for the dual space $Q^\ast$ in $\mathbb{R}^{(m+1)nl\ast}$, then we have $(\Psi_{su},\Phi_{su})=I_{(m+1)nl-c}$, the identity matrix in $\mathbb{R}^{((m+1)nl-c)\times((m+1)nl-c)}$. We denote by $C_Q$ the constant matrix such that $C\Phi_{su}=\Phi_{su}C_Q$; its spectrum coincides with $\sigma(C)\setminus \Lambda$.

With the notations above, the generalized eigenspace $\tilde{P}$ of $\tilde{C}$ associated with $\tilde \Lambda$ satisfies $\tilde{P}={\rm span}\tilde{\Phi}_c$,
where $\tilde{\Phi}_c={\rm diag}(\Phi_c,I_p)$. The left generalized eigenspace of $\tilde C$ associated with $\tilde\Lambda$ is spanned by $\tilde{\Psi}_c={\rm diag}(\Psi_c,I_p)$. Denoted by $\tilde{J}$ the $(c+p)\times(c+p)$ matrix subject to $\tilde{C}\tilde{\Phi}_c=\tilde{\Phi}_c\tilde{J}$, we have $\tilde{J}={\rm diag}(J,I_p)$ and its spectrum coincides with $\tilde\Lambda$. Clearly, $\sigma(\tilde C)\setminus\tilde\Lambda=\sigma(C)\setminus\Lambda $. If we denote $\tilde Q$ the complementary space of $\tilde P$ in the space $\mathbb{R}^{(m+1)nl+p}$, then $\tilde Q$ is spanned by $\tilde\Phi_{su}$, the generalized eigenspace of $\tilde C$ associated with $\sigma(\tilde C)\setminus \tilde\Lambda$. Obviously, $\tilde \Phi_{su}=(\Phi_{su}^T,0_p)^T$. We denote by $\tilde\Psi_{su}$ the left generalized eigenspace of $\tilde C$ associated with $\sigma(\tilde C)\setminus \tilde\Lambda$, then $\tilde \Psi_{su}=(\Psi_{su},0_p)$. Obviously, we have $(\tilde\Psi_{su},\tilde \Phi_{su})=I_{(m+1)nl-c}$.

According to the decomposition of the state space $\mathbb{R}^{(m+1)nl}=P\oplus Q$, $u$ could be represented as $u=\Phi_c x +\Phi_{su}y$ with $x\in\mathbb{R}^c$, $y\in \mathbb{R}^{(m+1)nl-c}$.




Therefore, projecting map (\ref{suspend}) onto $\tilde P$ and $\tilde Q$ respectively, we get
\begin{equation}\label{m10}
\begin{array}{rl}
\dbinom{x}{\alpha}&\hskip-3mm\mapsto\tilde J\dbinom{x}{\alpha}+\dbinom{\Psi_cF(\Phi_c x+\Phi_{su}y,\alpha)}{0_p},\\[3mm]
y&\hskip-3mm\mapsto C_Q y+\Psi_{su}F(\Phi_c x+\Phi_{su}y,\alpha).
\end{array}
\end{equation}
Noting that $\alpha$ is mapped into itself, the previous map exactly is
\begin{equation}\label{simp}
\begin{array}{l}
x\mapsto Jx+\Psi_cF(\Phi_c x+\Phi_{su}y,\alpha),\\
y\mapsto C_Q y+\Psi_{su}F(\Phi_c x+\Phi_{su}y,\alpha),
\end{array}
\end{equation}
where $\alpha$ should be considered as a variable. In addition we recall that the matrices $J$ and $C_Q$ are of the form of uptriangular.


The first part in map $(\ref{simp})$ has a fixed dimensionality, that is $c$, the dimensionality of the center manifold of (\ref{map}) associated with $\Lambda$. While the second part has a variable dimensionality, which tends to infinity as the parameter $\varepsilon$ tends to $0$. This makes it very hard to give explicitly the expressions for the (generalized) eigenvectors of (\ref{linmap}) associated with $\sigma(C)\setminus \Lambda$. Besides, the amount of the eigenvectors will also tend to infinity as the parameter $\varepsilon$ tends to $0$. It is impracticable to get the map restricted to the center manifold of (\ref{simp}) by the technique described in \cite{Wi90} directly. In the next we go deep into develop a method which allows getting the normal forms associated with $\Lambda$ on the center manifold avoid computing the eigenvectors $\Phi_{su}$ and $\Psi_{su}$.

Taylor expanding $F(u,\alpha)$ with respect to $(u,\alpha)$ leads to
\begin{equation}
F(u,\alpha)=\sum_{j\geq 2}\displaystyle\frac{1}{j!}F_j(u,\alpha), \ (u,\alpha)\in \mathbb{R}^{(m+1)nl+p}.
\end{equation}
Defining $f_j=(f_j^1,f_j^2)$ with
\begin{equation}\label{expand}
\begin{array}{rl}
f_j^1(x,y,\alpha)&=\Psi_cF_j(\Phi_cx+\Phi_{su}y,\alpha),\\
f_j^2(x,y,\alpha)&=\Psi_{su}F_j(\Phi_cx+\Phi_{su}y,\alpha),
\end{array}
\end{equation}
the map (\ref{suspend}) is equivalent to
\begin{equation}\label{seri}
\begin{array}{l}
x\mapsto J x+\sum\limits_{j\geq 2}\displaystyle\frac{1}{j!}f_j^1(x,y,\alpha),\\
y\mapsto C_Q y+\sum\limits_{j\geq 2}\displaystyle\frac{1}{j!}f_j^2(x,y,\alpha),
\end{array}
\end{equation}
where $x\in\mathbb{R}^c$ and $y\in\mathbb{R}^{(m+1)nl-c}$.

The normal forms are obtained by a recursive procedure, computing at each step the terms of order $j\geq 2$ in the normal form from the terms of the same order in the original map and the terms of lower orders already computed for the normal form in previous steps, through a transformation of variables
\begin{equation}\label{changecoor}
(x,y)=(\hat x,\hat y)+\displaystyle\frac{1}{j!}U_j(\hat x,\alpha),
\end{equation}
with $x,\ \hat x\in\mathbb{R}^c$, $y, \ \hat y\in\mathbb{R}^{(m+1)nl-c}$, and $U_j=(U_j^1, U_j^2)\in V_j^{c+p}(\mathbb{R}^c)\times V_j^{c+p}(\mathbb{R}^{(m+1)nl-c})$, where, for a normed space $X$, $V_j^{c+p}(X)$ denotes the linear space of homogeneous polynomials of degree $j$ in $c+p$ real variables, $(x,\alpha)=(x_1,x_2,\cdots,x_c,\alpha_1,\cdots,\alpha_p)$, and with coefficients in $X$,
\[
V_j^{c+p}(X)=\{\sum_{|(q,l)|=j} c_{(q,l)}x^q\alpha^l:(q,l)\in\mathbb{N}_0^{c+p}, c_{(q,l)}\in X\},
\]
$x^q\alpha^l=x_1^{q_1}\cdots x_{c}^{q_c}\alpha_1^{l_1}\cdots\alpha_p^{l_p}$ for $q=(q_1,\cdots,q_c)\in\mathbb{N}_0^c$ and $l=(l_1,\cdots,l_p)\in\mathbb{N}_0^p$, with the norm
\[
\left|\sum_{|(q,l)|=j}c_{(q,l)}x^q\alpha^l\right|=\sum_{|(q,l)|=j}|c_{(q,l)}|_X.
\]

We assume that after computing the normal form up to terms of order $j-1$ the map is
\begin{equation}\label{norm1}
\begin{array}{l}
x\mapsto Jx+\sum\limits_{i=2}^{j-1}\displaystyle\frac{1}{i!}g_i^1(x,y,\alpha)+\displaystyle\frac{1}{j!}\bar f_j^1(x,y,\alpha)+\cdots\\
y\mapsto C_Q y+\sum\limits_{i=2}^{j-1}\displaystyle\frac{1}{i!}g_i^2(x,y,\alpha)+\displaystyle\frac{1}{j!}\bar f_j^2(x,y,\alpha)+\cdots
\end{array}
\end{equation}
where
\[\begin{array}{l}
g_j^1(x,y,\alpha)=\bar f_j^1(x,y,\alpha)-[U_j^1(Jx,\alpha)-JU_j^1(x,\alpha)],\\[2mm]
g_j^2(x,y,\alpha)=\bar f_j^2(x,y,\alpha)-[U_j^2(Jx,\alpha)-C_QU_j^2(x,\alpha)].
\end{array}
\]
These formulas can be written for $g_j=(g_j^1,g_j^2)$ as
\begin{equation}\label{gf}
g_j=\bar f_j-M_jU_j,
\end{equation}
with $\bar f_j=(\bar f_j^1,\bar f_j^2)$ and $M_j$ defined below.

\begin{definition}
For $j\geq 2$, let $M_j$ denote the operator defined in $V_j^{c+p}(\mathbb{R}^c\times\mathbb{R}^{(m+1)nl-c})$, with values in the same space, by
\[
\begin{array}{rl}
M_j(p,h)&=(M_j^1p,M_j^2h)\\
(M_j^1p)(x,\alpha)&=p(Jx,\alpha)-Jp(x,\alpha)\\
(M_j^2h)(x,\alpha)&=h(Jx,\alpha)-C_Qh(x,\alpha),
\end{array}
\]
with domain $D(M_j)=V_j^{c+p}(\mathbb{R}^c)\times V_j^{c+p}(\mathbb{R}^{(m+1)nl-c})$.
\end{definition}

Then the problem we met is how to choose $U_j$ so that $g_j$ has a simple form.

We decompose the spaces $V_j^{c+p}(\mathbb{R}^c)$ and $
V_j^{c+p}(\mathbb{R}^{(m+1)nl-c})$ as
\[
\begin{array}{l}
V_j^{c+p}(\mathbb{R}^c)={\rm Im}(M_j^1)\oplus {\rm Im}(M_j^1)^c,\\
V_j^{c+p}(\mathbb{R}^c)=\ker(M_j^1)\oplus \ker(M_j^1)^c,
\end{array}
\]
and
\[\begin{array}{l}
V_j^{c+p}(\mathbb{R}^{(m+1)nl-c})={\rm Im}(M_j^2)\oplus{\rm Im}(M_j^2)^c,\\
V_j^{c+p}(\mathbb{R}^{(m+1)nl-c})=\ker(M_j^2)\oplus\ker(M_j^2)^c.
\end{array}\]
We denote the projections associated with the above decompositions of $V_j^{c+p}(\mathbb{R}^c)\times V_j^{c+p}(\mathbb{R}^{(m+1)nl-c})$ over ${\rm Im}(M_j^1)\times{\rm Im}(M_j^2)$ and over $\ker(M_j^1)^c\times\ker(M_j^2)^c$ by, respectively, $P_{I,j}=(P_{I,j}^1, P_{I,j}^2)$ and $P_{K,j}=(P_{K,j}^1,P_{K,j}^2)$. The complementary spaces in the above decompositions ${\rm Im}(M_j^i)^c$, ${\ker}(M_j^i)^c$ $(i=1,2)$, are not uniquely determined. As a consequence, normal forms are not unique, and depend on the choices of ${\rm Im}(M_j^i)^c$ $(i=1,2)$.

Let us now consider the right inverse of $M_j$ with range defined by the spaces complementary to the kernels of $M_j^i$ $(i=1,2)$, namely  $M_j^{-1}=((M_j^1)^{-1}, (M_j^2)^{-1})$ with $M_j^{-1}\circ P_{I,j}\circ M_j=P_{K,j}$.

Taking $y=0$ in formula (\ref{gf}), an adequate choice of $U_j$,
\begin{equation}
U_j(x,\alpha)=M_j^{-1}P_{I,j}\bar f_j(x,0,\alpha),
\end{equation}
allows taking away from $\bar f_j$ its component in the range of $M_j$, leading to
\begin{equation}
g_j(x,0,\alpha)=(I-P_{I,j})\bar f_j(x,0,\alpha).
\end{equation}

Therefore, the normal form for map (\ref{map}) relative to the invariant space $P$ and the projections $P_{I,j}$, $P_{K,j}$ $(j=2,3,\cdots)$ is the map in $\mathbb{R}^c\times\mathbb{R}^{(m+1)nl-c}$
\begin{equation}\label{normal}
\begin{array}{l}
x\mapsto Jx+\sum_{j\geq 2}\displaystyle\frac{1}{j!}g_j^1(x,y,\alpha)\\
y\mapsto C_Qy+\sum_{j\geq 2}\displaystyle\frac{1}{j!}g_j^2(x,y,\alpha).
\end{array}
\end{equation}

In the next, we show that the normal form relative to $P$ takes a more simple form of $c$-dimensional map on a locally invariant manifold for (\ref{map}) tangent to $P$ at zero.

Recall that the matrices $J$ and $C_Q$ are of Jordan forms, then by the Theorem 3.8 in \cite{ac88} we know that the normal form (\ref{normal}) can be chosen so that its nonlinear part contains only resonant monomials.
\begin{definition}
We say that the map (\ref{suspend}) satisfies the nonresonance conditions relative to $\tilde \Lambda\in\sigma(\tilde C)$, if
\[
\bar\lambda^q\neq\mu,\ \ \ \mbox{for all } \mu \in\sigma(\tilde C)\setminus\tilde\Lambda, \ q\in\mathbb{N}_0^{c+p},\ |q|\geq 2,
\]
where $\bar \lambda=(\lambda_1,\cdots,\lambda_{c+p})$ and $\lambda_1,\cdots,\lambda_{c+p}$ are the elements of $\tilde\Lambda$, each one of them appearing as many times as its multiplicity as an eigenvalue of the matrix $\tilde C$.
\end{definition}

Clearly, the fact that the map (\ref{suspend}) satisfies the nonresonance conditions relative to $\tilde \Lambda\in\sigma(\tilde C)$ is equivalent to the map (\ref{map}) satisfies the nonresonance condtions relative to $\Lambda\in\sigma( C)$ since the difference between $\tilde \Lambda$ and $\Lambda$ are the eigenvalues introduced by regarding $\alpha$ as a variable, that is $\tilde\Lambda\setminus\Lambda=\{\underbrace{1,\cdots,1}_{p\mbox{ }times}\}$.

\begin{theorem}\label{mainth}
Consider the map (\ref{map}) and let $P$ be the invariant subspace of $C$ associated with a nonempty finite set $\Lambda$ of eigenvalues. For the decomposition of $\mathbb{R}^{(m+1)nl}=\mathbb{R}^c\oplus \mathbb{R}^{(m+1)nl-c}$, we get $u=\Phi_c x+\Phi_{su}y$, where $x\in\mathbb{R}^c$ and $y\in\mathbb{R}^{(m+1)nl-c}$, $\Phi_c$ and $\Phi_{su}$ are formed by the generalized eigenvectors of $C$ associated with $\Lambda$ and $\sigma(C)\setminus \Lambda$, respectively. If the nonresonance conditions relative to $\Lambda$ are satisfied, then there exists a formal change of variables $(\bar x,\bar y)\mapsto (x,y)$ of the form $x=\bar x+p(\bar x,\alpha), y=\bar y+h(\bar x,\alpha)$, such that the map (\ref{map}) is equivalent to the following map in $\mathbb{R}^c\times\mathbb{R}^{(m+1)nl-c}$
\begin{equation}
\begin{array}{l}
\bar x\mapsto J\bar x+\sum_{j\geq 2}\displaystyle\frac{1}{j!}g_j^1(\bar x,\bar y,\alpha)\\
\bar y\mapsto C_Q \bar y+\sum_{j\geq 2}\displaystyle\frac{1}{j!}g_j^2(\bar x,\bar y,\alpha),
\end{array}
\end{equation}
where $g_j^1,g_j^2$ are computed as (\ref{gf}), with $g_j^2(\bar x,0, \alpha)=0$ for all $j\geq 2$. This map is in normal form relative to $P$. If there exists a locally invariant manifold for the map (\ref{map}) tangent to $P$ at zero, then it satisfies $\bar y=0$ and the map on it is given by the $c$-dimensional map
\begin{equation}
\bar x\mapsto J\bar x+\sum_{j\geq 2}\displaystyle\frac{1}{j!}g_j^1(\bar x,0, \alpha),
\end{equation}
which is in normal form for maps.
\end{theorem}

\section{Takens-Bogdanov bifurcations of parameterized delay differential equations}

To show the bifurcation structure near Takens-Bogdanov point could be preserved by numerical discretization, we recall the results on the bifurcation structure near Takens-Bogdanov point in this section, cf. \cite{xu08}.

We consider the DDEs of the type
\begin{equation}\label{ep1}
\dot{z}(t)=f(z(t),z(t-1),\alpha),
\end{equation}
where $z\in \mathbb{R}^n$, $\alpha\in\mathbb{R}^2$ is a bifurcation parameter, $f(z_1,z_2,\alpha)$ is a $C^r(r\geq2)$
smooth function from
$\mathbb{R}^n\times\mathbb{R}^n\times\mathbb{R}^2$ to
$\mathbb{R}^n$ with
\begin{equation}\label{e2.2.3}
f(0,0,\alpha)=0,\dfrac{\partial{f}}{\partial{z_1}}(0,0,\alpha)=0,\dfrac{\partial{f}}{\partial{z_2}}(0,0,\alpha)=0,
\ \ \forall \alpha\in\mathbb{R}^2.
\end{equation}
Denote $A=\frac{\partial f}{\partial z_1}(0,0,0)$ and $B=\frac{\partial f}{\partial z_2}(0,0,0)$, then we can rewrite equation
(\ref{ep1}) as
\begin{equation}\label{e2.2.5}
\dot{z}(t)=A z(t)+B z(t-1)+F(z(t),z(t-1),\alpha),
\end{equation}
which could be linearized at $(z,\alpha)=(0,0)$ as
\begin{equation}\label{elin}
\dot{z}(t)=A z(t)+B z(t-1).
\end{equation}
The characteristic equation of (\ref{elin}) is given by
\begin{equation}\label{cha}\det(\mu I-A-Be^{-\mu})=0.\end{equation}
We assume $(z,\alpha)=(0,0)$ is a Takens-Bogdanov point of (\ref{ep1}), i.e.,
\begin{description}
\item[A]$0$ is a root of (\ref{cha}) with algebraic multiplicity 2 and geometric multiplicity 1, and the other roots exhibit nonzero  real parts.
\end{description}


\begin{lemma}{\rm \cite{xu08}}\label{l2.3.1}There exist $\phi_1^0\in\mathbb{R}^n\backslash\{0\}$,
 $\phi_2^0\in\mathbb{R}^n$, and
 $\psi_2^0\in\mathbb{R}^{n\ast}\backslash\{0\}$, $\psi_1^0\in\mathbb{R}^{n\ast}$, they satisfy the following equations
 \begin{equation}\label{e2.3.7}
\begin{array}{l}
\begin{array}{ll}
(1)\ (A+B)\phi_1^0=0,&\qquad\qquad(2)\ (A+B)\phi_2^0=(B+I)\phi_1^0,\\
(3)\ \psi_2^0(A+B)=0,&\qquad\qquad(4)\ \psi_1^0(A+B)=\psi_2^0(B+I),\\
\end{array}\\
\begin{array}{l}
(5)\
\psi_2^0\phi_2^0-\frac{1}{2}\psi_2^0B\phi_1^0+\psi_2^0B\phi_2^0=1,\\[1mm]
(6)\
\psi_1^0\phi_2^0-\frac{1}{2}\psi_1^0B\phi_1^0+\psi_1^0B\phi_2^0
+\frac{1}{6}\psi_2^0B\phi_1^0-\frac{1}{2}\psi_2^0B\phi_2^0=0.
\end{array}
\end{array}
\end{equation}
\end{lemma}
Taylor expanding  $F(z(t),z(t-1),\alpha)$ with respect to
$z(t)$, $z(t-1)$ and $\alpha$ we obtain
\begin{equation}\label{e2.3.4}
{F}(z(t),z(t-1),\alpha)=\sum\limits_{j\geq2}\dfrac{1}{j!}
F_j(z(t),z(t-1),\alpha),
\end{equation}
where the first term $(j=2)$ can be expressed in the form
\begin{equation}\label{ch2fe3.16}
\begin{array}{l}
\dfrac{1}{2}F_2(z(t),z(t-1),\alpha)\\
=A_1\alpha_1z(t)+A_2\alpha_2z(t)+B_1\alpha_1z(t-1)+B_2\alpha_2z(t-1)\\
\quad+\sum_{i=1}^nE_iz_i(t)z(t-1)+\sum_{i=1}^nF_iz_i(t)z(t)\\
\quad+\sum_{i=1}^nG_iz_i(t-1)z(t-1)
\end{array}
\end{equation}
with $A_i,B_i(i=1,2),E_i,F_i,G_i(i=1,2,\cdots,n)$ coefficient
matrices, and there is no terms of $O(\alpha^2)$ in $
F_2(z(t),z(t-1),\alpha)$ since ${F}(0,0,\alpha)=0,\forall
\alpha\in\mathbb{R}^2$.

Denote
\begin{equation}\label{ab}
\begin{array}{rl}
a=&\psi_2^0\sum_{i=1}^{n}(E_i+F_i+G_i)\phi_1^0\phi_{1i}^0,\\[3mm]
b=&2\psi_1^0\sum_{i=1}^n(E_i+F_i+G_i)\phi_1^0\phi_{1i}^0\\
&+
\psi_2^0\{\sum_{i=1}^n(E_i+F_i+G_i)(\phi_2^0\phi_{1i}^0+
\phi_1^0\phi_{2i}^0)-\sum_{i=1}^n(E_i+2G_i)\phi_1^0\phi_{1i}^0\},
\end{array}
\end{equation}
and
\begin{equation}\label{k1k2}
\dbinom{\kappa_1}{\kappa_2}= \Pi\dbinom{\alpha_1}{\alpha_2}
\end{equation}
with
$$ \Pi={\left(\begin{array}{ll}
\ \ \psi_2^0(A_1+B_1)\phi_1^0&\ \ \psi_2^0(A_2+B_2)\phi_1^0\\[2mm]
\begin{array}{l}\{\psi_1^0(A_1+B_1)\phi_1^0\\+\psi_2^0((A_1+B_1)\phi_2^0-B_1\phi_1^0)\}\end{array}&
\begin{array}{l}\{\psi_1^0(A_2+B_2)\phi_1^0\\+\psi_2^0((A_2+B_2)\phi_2^0-B_2\phi_1^0)\}\end{array}
\end{array}\right)}.
$$
For the bifurcation structure near the Takens-Bogdanov point of DDE (\ref{ep1}), we have the following theorem.

\begin{theorem}{\rm \cite{xu08}}\label{t2.4.1}{Assume the assumption {\bf A} holds, $\det \Pi\neq0$
and $a\cdot b\neq0$. Then there exists a constant $\kappa_1^0>0$, such
that when $0<\kappa_1(\alpha_1,\alpha_2)<\kappa_1^0$, in the
parameter plane $(\alpha_1,\alpha_2)$ near the origin there exist
two curves:
${l}_h$ and ${l}_\infty$\\
\indent 1.\ \  the curve ${l}_h$, which has the
following local representation:
$$
{l}_h=\{(\alpha_1,\alpha_2):\
\kappa_2(\alpha_1,\alpha_2)-\dfrac{b}{a}\kappa_1(\alpha_1,\alpha_2)+h.o.t.=0,\
\kappa_1(\alpha_1,\alpha_2)>0\},
$$
is a Hopf point branch of the DDE
(\ref{ep1}), where $h.o.t.=o(|(\alpha_1,\alpha_2)|)$, i.e.
$ l_h$ consists of
Hopf bifurcation points of (\ref{ep1});\\
\indent 2.\ \  the curve ${l}_\infty$, which has the
following local representation:
\begin{equation}\label{e2.4.6}
{l}_\infty=\{(\alpha_1,\alpha_2):\
h(\alpha_1,\alpha_2)+h.o.t.=0,\ \kappa_1(\alpha_1,\alpha_2)>0\},
\end{equation}
is a homoclinic branch of the DDE
(\ref{ep1}), where
$h(\alpha_1,\alpha_2)=\kappa_2(\alpha_1,\alpha_2)
-\mu(\sqrt{\kappa_1(\alpha_1,\alpha_2)})\linebreak[3]
\kappa_1(\alpha_1,\alpha_2)$, $\mu(\cdot)$ is a continuously differentiable function with $\mu(0)=\frac{6}{7}ba^{-1}$ and $h.o.t.=o(|(\alpha_1,\alpha_2)|)$. In other
words, equation (\ref{ep1}) has a unique homoclinic orbit
connecting the origin for each
$(\alpha_1,\alpha_2)\linebreak[0]\in{l}_\infty$. }
\end{theorem}

\section{ Preservation of Takens-Bogdanov bifurcations by Euler discretization}

The forward Euler scheme for solving (\ref{ep1}) with step size $\varepsilon=\frac{1}{m}, m\in \mathbb{Z}^+$ is given by
\begin{equation}\label{add}
z_{k+1}=z_k+\varepsilon f(z_k,z_{k-m},\alpha),
\end{equation}
which can be reformulated to
\begin{equation}\label{Euler}
z_{k+1}=z_k+\varepsilon A z_k+\varepsilon B z_{k-m}+\varepsilon F(z_k, z_{k-m},\alpha).
\end{equation}

Denoting $u_k=(z_k^T,z_{k-1}^T,\dots,z_{k-m}^T)^T\in\mathbb{R}^{(m+1)n}$, (\ref{Euler}) can be rewritten as
\begin{equation}\label{Euler2}
u_{k+1}=C u_k+H(u_k,\alpha),
\end{equation}
where
$$
C=\left(\begin{array}{ccccc}
I+\varepsilon A& 0 &\dots& 0 & \varepsilon B\\
I& 0 &\dots& 0 & 0\\
0&\ddots&\ddots&\vdots&\vdots\\
\vdots&\ddots&\ddots&0&\vdots\\
0&\cdots&0 & I & 0\end{array}
\right)
$$
and
$$
H(u_k,\alpha)=(\varepsilon F(z_k,z_{k-m},\alpha)^T,0, \cdots, 0)^T.
$$

Linearizing the map (\ref{Euler2}) at $(u,\alpha)=(0,0)$ we obtain
\begin{equation}\label{linearmap}
u_{k+1}=C u_k,
\end{equation}
which has the characteristic equation
\begin{equation}\label{chara}
\det (\lambda I-C)=0.
\end{equation}
Noting (\ref{Euler2}) is a map of the type of (\ref{map}) as $l=1$, the techniques developed in Section 2 allow us to get the normal forms on the low dimensional center manifold without computing the eigenvectors associated with the eigenvalues of the modulus other than $1$. While the map (\ref{map}) as $l>1$ could correspond to more general numerical schemes, eg. the general one step method.

Under the assumption that (\ref{ep1}) undergoes a Takens-Bogdanov bifurcation at $(z,\alpha)=(0,0)$, we first show the Takens-Bogdanov point $(z,\alpha)=0,0)$ of (\ref{ep1}) is inherited without any shift by the Euler method (\ref{Euler}) and turns into a 1:1 resonance point.
\begin{theorem}\label{t41}
Assume the assumption {\bf A} holds. Then (\ref{add}) undergoes a 1:1 resonance at $(z,\alpha)=(0,0)$.
\end{theorem}
\proofbegin
We only need to show the double zero eigenvalue of (\ref{elin}) turns into a double unit multiplier, $\lambda_{1,2}=1$ at $\alpha=0$, and no other eigenvalue of (\ref{linearmap}) has modulus 1.

Assume
\begin{equation}\begin{array}{rlrl}(C-I)\phi_{1}&=0,&(C-I)\phi_{2}&=\phi_{1},\\
\psi_{2}(C-I)&=0,&\psi_{1}(C-I)&=\psi_{2},\end{array}\end{equation}
and require
$$
\psi_1\cdot\phi_1=\psi_2\cdot\phi_2=1.
$$
We can check that the following choice of  $\phi_1,\phi_2,\psi_2,\psi_1$
\begin{equation}\begin{array}{l}
\phi_1=\varepsilon(\phi_1^{0T},\cdots,\phi_1^{0T})^T,\\
\phi_2=\varepsilon(m\phi_2^{0T},m\phi_2^{0T}-\phi_1^{0T},\cdots,m\phi_2^{0T}-m\phi_1^{0T})^T,\\
\psi_2=\dfrac{1}{1-\frac{1}{2m}\psi_2^0B\phi_1^0}(\psi_2^0,\varepsilon\psi_2^0B,\cdots,\varepsilon\psi_2^0B),\\
\psi_1=\dfrac{1}{1-\frac{1}{2m}\psi_2^0B\phi_1^0}(m\psi_1^0,\psi_1^0B-\psi_2^0B,\psi_1^0B-(m-1)\varepsilon\psi_2^0B,\cdots,\psi_1^0B-\varepsilon\psi_2^0B),
\end{array}\end{equation}
fulfill the requirements above. Besides, the Fredholm alternative Theorem implies
$$
\psi_2\cdot\phi_1=\psi_1\cdot\phi_2=0.
$$
These show that 1 is an eigenvalue of  (\ref{linearmap}) with algebraic multiplicity 2 and geometric multiplicity 1.

Next, we show no other eigenvalue of (\ref{linearmap}) has modulus 1.

Denote $d(\mu)=e^{\mu}(\mu I-A)-B $ and $D(\mu,\varepsilon)=e^\mu(g(\mu\varepsilon)\mu I-A)-B$, where $g(x)=\frac{e^x-1}{x}$. Then, see \cite{xu10}, we have $\det d(\mu)=0$ is equivalent to (\ref{cha}), and if $\det D(\mu_\varepsilon,\varepsilon)=0$, then $\lambda=e^{\mu_\varepsilon \varepsilon}$ solves  (\ref{chara}). Besides, we have $\lim\limits_{\varepsilon\rightarrow 0}D(\mu,\varepsilon)=d(\mu)$. Hence for any series $\mu_\varepsilon$ of solutions of $\det D(\mu,\varepsilon)=0$, there exists a solution $\mu_0$ of $\det d(\mu)=0$ such that $\lim\limits_{\varepsilon\rightarrow 0}\mu_\varepsilon=\mu_0$. Therefore, if other than the eigenvalue 1 there exists any eigenvalue of (\ref{chara}), denoted by $\lambda_\varepsilon$, s.t. $|\lambda_\varepsilon|=1$, then $\mu_\varepsilon=\frac{1}{\varepsilon}\ln \lambda_\varepsilon$ solves $\det D(\mu,\varepsilon)=0$. Since $\mu_\varepsilon$ is pure imaginary, we have $\mu_0$, as the limit of $\mu_\varepsilon$ as $\varepsilon\rightarrow0$, is either pure imaginary or 0 and solves (\ref{cha}). This contradicts to the assumption {\bf A}.
$\Box$


Let $\Phi_c=(\phi_1,\phi_2)$, $\Psi_c=(\psi_1,\psi_2)$, and $\Lambda=\{1,1\}$, then we know the invariant space $P$  of (\ref{Euler2}) associated with $\Lambda$ is spanned by $\Phi_c$, the dual invariant space $P^\ast$ of (\ref{Euler2}) associated with $\Lambda$ is spanned by $\Psi_c$. They satisfy $(\Psi_c,\Phi_c)=I_2$, $C\Phi_c=\Phi_c J$, and $\Psi_c C=J\Psi_c$ with $J=\left(\begin{array}{cc}1&1\\ 0&1\end{array}\right)$.

We then consider the Taylor expansion of $H(u_k,\alpha)$ with respect to $(u_k,\alpha)$. By (\ref{e2.3.4}), we have
\begin{equation}
H(u_k,\alpha)=\varepsilon({F}(z_k,z_{k-m},\alpha)^T,0,\cdots,0)^T=(\sum\limits_{j\geq2}\dfrac{\varepsilon}{j!}
F_j(z_k,z_{k-m},\alpha)^T,0,\cdots,0)^T,
\end{equation}
where
\begin{equation}\label{taylor2}
\begin{array}{l}
\dfrac{1}{2}F_2(z_k,z_{k-m},\alpha)\\
=A_1\alpha_1z_k+A_2\alpha_2z_k+B_1\alpha_1z_{k-m}+B_2\alpha_2z_{k-m}\\
\quad+\sum_{i=1}^nE_iz_k^iz_{k-m}+\sum_{i=1}^nF_iz_k^iz_k\\
\quad+\sum_{i=1}^nG_iz_{k-m}^iz_{k-m}.
\end{array}
\end{equation}
Evidently, the canonical basis of $V_2^4(\mathbb{R}^2)$ is
composed by the elements
$$
\begin{array}{l}
\dbinom{x_1^2}{0},\ \dbinom{x_2^2}{0},\ \dbinom{\alpha_1^2}{0},\
\dbinom{\alpha_2^2}{0},\ \dbinom{x_1x_2}{0},\
\dbinom{x_1\alpha_1}{0},\ \dbinom{x_1\alpha_2}{0},
\dbinom{x_2\alpha_1}{0},\\[4mm]
\dbinom{x_2\alpha_2}{0},\ \dbinom{\alpha_1\alpha_2}{0},
\dbinom{0}{x_1^2},\ \dbinom{0}{x_2^2},\ \dbinom{0}{\alpha_1^2},\
\dbinom{0}{\alpha_2^2},\ \dbinom{0}{x_1x_2},\
\dbinom{0}{x_1\alpha_1},\\[4mm]
 \dbinom{0}{x_1\alpha_2},\
\dbinom{0}{x_2\alpha_1},\ \dbinom{0}{x_2\alpha_2},\
\dbinom{0}{\alpha_1\alpha_2},\
\end{array}
$$
and the images of these elements under $M_2^1$ are, respectively
$$
\begin{array}{l}
\dbinom{2x_1x_2+x_2^2}{0},\ \dbinom{x_2^2}{0},\ \dbinom{0}{0},\
\dbinom{0}{0},\ \dbinom{x_2^2}{0},\ \dbinom{x_2\alpha_1}{0},\
\dbinom{x_2\alpha_2}{0}, \dbinom{0}{0},\\[4mm]
\dbinom{0}{0},\ \dbinom{0}{0},\ \dbinom{-x_1^2}{2x_1x_2+x_2^2},\
\dbinom{-x_2^2}{0}, \dbinom{-\alpha_1^2}{0},\
\dbinom{-\alpha_2^2}{0},\ \dbinom{-x_1x_2}{x_2^2},\
\dbinom{-x_1\alpha_1}{x_2\alpha_1},\\[4mm]
\dbinom{-x_1\alpha_2}{x_2\alpha_2},\ \dbinom{-x_2\alpha_1}{0},\
\dbinom{-x_2\alpha_2}{0},\ \dbinom{-\alpha_1\alpha_2}{0}.\
\end{array}
$$
Therefore, a basis of ${\rm Im}(M_2^1)^c$ can be taken as the
set composed by the elements
$$
\begin{array}{l}
\dbinom{0}{x_1^2},\ \dbinom{0}{\alpha_1^2},\
\dbinom{0}{\alpha_2^2},\ \dbinom{0}{x_1x_2},\
\dbinom{0}{x_1\alpha_1},\ \dbinom{0}{x_1\alpha_2},\
\dbinom{0}{x_2\alpha_1},\ \dbinom{0}{x_2\alpha_2},\
\dbinom{0}{\alpha_1\alpha_2}.
\end{array}
$$

Denoting by $\Phi_c^0$ the first $n$ rows of $\Phi_c$, $\Phi_c^m$ the last $n$ rows of $\Phi_c$, we have
$H(\Phi_cx,\alpha)=(\varepsilon F(\Phi_c^0x,\Phi_c^mx,\alpha)^T,0,\cdots,0)^T.$

Denoting by $\phi_{ji}$ the $i$-th element of $\phi_j$, we have
from (\ref{ch2fe3.16})
\begin{eqnarray*}
&&\dfrac{1}{2}{F}_2((\varepsilon\phi_1^0,\phi_2^0)x,(\varepsilon\phi_1^0,\phi_2^0-\phi_1^0)x,\alpha)\\
&=&A_{1}\alpha_1(\varepsilon\phi_1^0,\phi_2^0)(x_1,x_2)^T+
A_{2}\alpha_2(\varepsilon\phi_1^0,\phi_2^0)(x_1,x_2)^T\\
&&+B_{1}\alpha_1(\varepsilon\phi_1^0,\phi_2^0-\phi_1^0)(x_1,x_2)^T+
B_{2}\alpha_2(\varepsilon\phi_1^0,\phi_2^0-\phi_1^0)(x_1,x_2)^T\\
&&+\sum_{i=1}^nE_i(\varepsilon\phi_{1i}^0,\phi_{2i}^0)(x_1,x_2)^T(\varepsilon\phi_1^0,\phi_2^0-\phi_1^0)(x_1,x_2)^T\\
&&+\sum_{i=1}^n
F_i(\varepsilon\phi_{1i}^0,\phi_{2i}^0)(x_1,x_2)^T(\varepsilon\phi_1^0,\phi_2^0)(x_1,x_2)^T\\
&&+\sum_{i=1}^n
G_i(\varepsilon\phi_{1i}^0,\phi_{2i}^0-\phi_{1i}^0)(x_1,x_2)^T(\varepsilon\phi_1^0,\phi_2^0-\phi_1^0)(x_1,x_2)^T\\
&=&\varepsilon(A_1+B_1)\phi_1^0\alpha_1x_1+\varepsilon(A_2+B_2)\phi_1^0\alpha_2x_1\\
&&+((A_1+B_1)\phi_2^0-B_1\phi_1^0)\alpha_1x_2+((A_2+B_2)\phi_2^0-B_2\phi_1^0)\alpha_2x_2\\
&&+\varepsilon^2\sum_{i=1}^n(E_i+F_i+G_i)\phi_{1i}^0\phi_1^0)x_1^2\\
&&+\varepsilon\sum_{i=1}^n\{(E_i+F_i+G_i)(\phi_{1i}^0\phi_{2}^0+\phi_{2i}^0\phi_1^0)-(E_i+2G_i)\phi_{1i}^0\phi_1^0\}x_1x_2\\
&&+\sum_{i=1}^n(E_i+F_i+G_i)\phi_{2i}^0\phi_{2}^0-(E_i+G_i)\phi_{2i}^0\phi_{1}^0-G_i\phi_{1i}^0(\phi_2^0-\phi_1^0))x_2^2.\\
\end{eqnarray*}

Base on the expansion above, and the canonical basis of $V_2^4(\mathbb{R}^2)$, ${\rm Im}(M_2^1)$ and ${\rm Im}(M_2^1)^c$, noting that $f_2^1(x,0,\alpha)=\Psi_cH_2(\Phi_cx,\alpha)=\dfrac{1}{2}(\frac{1}{1-\frac{1}{2m}\psi_2^0B\phi_1^0}\psi_1^0,\frac{\varepsilon}{1-\frac{1}{2m}\psi_2^0B\phi_1^0}\psi_2^0)\linebreak[3]
{F}_2((\varepsilon\phi_1^0,\phi_2^0)x,(\varepsilon\phi_1^0,\phi_2^0-\phi_1^0)x,\alpha)$, we can compute the function $g_2^1(x,0,\alpha)=(I-P_{I,2}^1)f_2^1(x,0,\alpha)$. By Theorem \ref{mainth}, we have the following results.

\begin{theorem}\label{t3.3}
Assume the requirements in Theorem \ref{t2.4.1} are fulfilled. Then the numerical scheme (\ref{add}) will undergo a 1:1 resonance at $(z,\alpha)=(0,0)$. In addition, the numerical scheme (\ref{add}) could be reduced to a 2 dimensional map on the center manifold at $(z,\alpha)=(0,0)$ as follows
\begin{equation}\label{reduced}
\begin{array}{l}
x_1\mapsto x_1+x_2,\\
x_2\mapsto x_2+\kappa_1^\varepsilon x_1+\kappa_2^\varepsilon x_2+a^\varepsilon x_1^2+b^\varepsilon x_1x_2+h.o.t.,
\end{array}
\end{equation}
where
$$\begin{array}{rlrl}
\kappa_1^\varepsilon&=\frac{\varepsilon^2}{1-\frac{\varepsilon}{2}\psi_2^0 B\phi_1^0}\kappa_1,&
\kappa_2^\varepsilon&=\frac{\varepsilon}{1-\frac{\varepsilon}{2}\psi_2^0 B\phi_1^0}\kappa_2,\\[3mm]
a^\varepsilon&=\frac{\varepsilon^3}{1-\frac{\varepsilon}{2}\psi_2^0 B\phi_1^0}a,&
b^\varepsilon&=\frac{\varepsilon^2}{1-\frac{\varepsilon}{2}\psi_2^0 B\phi_1^0}b,
\end{array}
$$
with $\kappa_1,\kappa_2$ and $a,b$ defined in (\ref{k1k2}) and (\ref{ab}), respectively.

\end{theorem}

It is known that for the reduced map (\ref{reduced}), if $a^\varepsilon\cdot b^\varepsilon\neq 0$ (equivalently $a\cdot b\neq 0$), the local bifurcation structure near $(x,\alpha)=(0,0)$ is determined by the linear and quadratic terms, and not the terms of order higher. Hence  we turn to investigate the local bifurcation structures of the map
\begin{equation}\label{reducedsimp}
\begin{array}{l}
x_1\mapsto x_1+x_2,\\
x_2\mapsto x_2+\kappa_1^\varepsilon x_1+\kappa_2^\varepsilon x_2+a^\varepsilon x_1^2+b^\varepsilon x_1x_2.
\end{array}
\end{equation}

\begin{lemma}\label{lemma1}
{Let $\lambda_\varepsilon^\pm(\kappa_1^\varepsilon,\kappa_2^\varepsilon)$ be the eigenvalues of
the Jacobian of (\ref{reducedsimp}) at $(-\frac{\kappa_1^\varepsilon}{a^\varepsilon},0)$. Then, when $a^\varepsilon\cdot b^\varepsilon\neq0$ and $0<\kappa_1^\varepsilon<2$\\
$$|\lambda_\varepsilon^\pm(\kappa_1^\varepsilon,\kappa_2^\varepsilon)|=\sqrt{1+\kappa_2^\varepsilon-\frac{b^\varepsilon}{a^\varepsilon}\kappa_1^\varepsilon+\kappa_1^\varepsilon}.$$
Hence we conclude that each point on the line segment
$\tilde l_h^\varepsilon=\{(\kappa_1^\varepsilon,\kappa_2^\varepsilon):\
\kappa_2^\varepsilon=\dfrac{b^\varepsilon}{a^\varepsilon}\kappa_1^\varepsilon-\kappa_1^\varepsilon,0<\kappa_1^\varepsilon<2 \}$ in the parameter plane
$(\kappa_1^\varepsilon$,$\kappa_2^\varepsilon)$ is a Neimark-Sacker bifurcation (also known as the Hopf bifurcation for map) point of the map (\ref{reducedsimp}).}
\end{lemma}
\proofbegin Evidently the map (\ref{reducedsimp}) has two fixed points, $(-\frac{\kappa_1^\varepsilon}{a^\varepsilon},0)$ and $(0,0)$. Direct computations show the eigenvalues of the Jacabian at $(-\frac{\kappa_1^\varepsilon}{a^\varepsilon},0)$ reads
$$ \lambda_\varepsilon^\pm(\kappa_1^\varepsilon,\kappa_2^\varepsilon)=\frac{1}{2}(2+\kappa_2^\varepsilon-\frac{b^\varepsilon}{a^\varepsilon}\kappa_1^\varepsilon)
\pm\sqrt{(2+\kappa_2^\varepsilon-\frac{b^\varepsilon}{a^\varepsilon}\kappa_1^\varepsilon)^2-
4(1+\kappa_2^\varepsilon-\frac{b^\varepsilon}{a^\varepsilon}\kappa_1^\varepsilon+\kappa_1^\varepsilon)},$$
which have a modulus given by
$$
|\lambda_\varepsilon^\pm(\kappa_1^\varepsilon,\kappa_2^\varepsilon)|=\sqrt{1+\kappa_2^\varepsilon-\dfrac{b^\varepsilon}{a^\varepsilon}\kappa_1^\varepsilon+\kappa_1^\varepsilon}.
$$
Therefore, when $(\kappa_1^\varepsilon,\kappa_2^\varepsilon)$ changes from one side of $\tilde l_h^\varepsilon$ to
the other, in the parameter plane the eigenvalues
$\lambda_\varepsilon^\pm(\kappa_1^\varepsilon,\kappa_2^\varepsilon)$ will cross the unit circle from
outside to inside ($\kappa_2^\varepsilon-\dfrac{b^\varepsilon}{a^\varepsilon}\kappa_1^\varepsilon+\kappa_1^\varepsilon>0$) or from inside to outside ($\kappa_2^\varepsilon-\dfrac{b^\varepsilon}{a^\varepsilon}\kappa_1^\varepsilon+\kappa_1^\varepsilon<0$).
Let $\theta=\arctan\frac{\sqrt{4\kappa_1^\varepsilon+(\kappa_1^\varepsilon)^2}}{2-\kappa_1^\varepsilon}$, obviously we have when $0<\kappa_1^\varepsilon<2$
$$e^{ik\theta}\neq 1,\ \ \ \mbox{for}\ k=1,2,3,4.$$
 Applying of the Neimark-Sacker bifurcation Theorem implies
$(x_1,x_2,\kappa_1^\varepsilon,\kappa_2^\varepsilon)=(-\dfrac{\kappa_1^\varepsilon}{a^\varepsilon},0,\linebreak[3]\kappa_1^\varepsilon,\dfrac{b^\varepsilon}{a^\varepsilon}\kappa_1^\varepsilon-\kappa_1^\varepsilon)$
is a Neimark-Sacker bifurcation point of the map (\ref{reducedsimp}) as
$0<\kappa_1^\varepsilon<2$. In fact, there is a neighborhood of $(-\frac{\kappa_1^\varepsilon}{a^\varepsilon},0)$ in which a unique closed invariant curve bifurcates from $(-\frac{\kappa_1^\varepsilon}{a^\varepsilon},0)$ when $\kappa_2^\varepsilon-\dfrac{b^\varepsilon}{a^\varepsilon}\kappa_1^\varepsilon+\kappa_1^\varepsilon$ changes signs. In other words, the line segment $\tilde l_h^\varepsilon$ in
the parameter plane $(\kappa_1^\varepsilon,\kappa_2^\varepsilon)$ is a Neimark-Sacker point branch
of the map (\ref{reducedsimp}).
$\Box$

\begin{remark} Noting the expression of $\kappa_1^\varepsilon\approx \varepsilon^2\kappa_1$, the statement of $\kappa_1^\varepsilon<2$ requires the step size $\varepsilon$ of the numerical scheme (\ref{add}) should be taken nicely small to reproduce the Hopf bifurcations of (\ref{ep1}).
\end{remark}

In the next we consider the homoclinic curves bifurcates from the fixed point $(x,\alpha)=(0,0)$.

Applying the transformation of $\bar x_1=a^\varepsilon(x_1+\frac{\kappa_1^\varepsilon}{2a^\varepsilon}), \bar x_2=a^\varepsilon x_2$ to (\ref{reducedsimp}) leads to another typical normal form for 1:1 resonance, that is
\begin{equation}
\begin{array}{l}
\bar x_1\mapsto\bar x_1+\bar x_2,\\
\bar x_2\mapsto\bar x_2-\frac{(\kappa_1^\varepsilon)^2}{4}+(\kappa_2^\varepsilon -\frac{b^\varepsilon\kappa_1^\varepsilon}{2a^\varepsilon}) \bar x_2+\bar x_1^2+\frac{b^\varepsilon}{a^\varepsilon} \bar x_1\bar x_2.
\end{array}
\end{equation}
Applying the results in \cite{Broer} to the former map shows that  the homoclinic curves bifurcated from the fixed points $(x,\alpha)=(0,0)$ is depicted by
\begin{equation}
\kappa_2^\varepsilon=\frac{6}{7}\frac{b^\varepsilon}{a^\varepsilon}\kappa_1^\varepsilon-\frac{5}{7}\kappa_1^\varepsilon+O((\kappa_1^\varepsilon)^\frac{3}{2})
\end{equation}
as $0<\kappa_1^\varepsilon<\kappa_\varepsilon^0$, where $\kappa_\varepsilon^0$ is some positive constant.

\begin{remark}
In fact, there exist two curves $\alpha_2^+(\alpha_1)$ and $\alpha_2^-(\alpha_1)$ respectively corresponding to the first and the last homoclinic tangency, they are exponentially close to one-another. If $\alpha$ is located in the region confined by these two curves, the map (\ref{reducedsimp}) possesses transverse homoclinic trajectories, see \cite{Gn09} for detail. But this is not our goal in this paper.
\end{remark}

Noting that the map (\ref{reducedsimp}) is locally topologically equivalent near the origin to (\ref{reduced}), we have the following bifurcation results based on Lemma \ref{lemma1} and the discussions above.

\begin{theorem}\label{t44}
{Assume $a^\varepsilon\cdot b^\varepsilon\neq0$ ($a^\varepsilon,b^\varepsilon\in\mathbb{R}$ defined in Theorem
\ref{t3.3}). Then there exists a constant $\kappa_{1\varepsilon}^0=\min\{2,\kappa_\varepsilon^0\}>0$, such
that when $0<\kappa_1^\varepsilon(\alpha_1,\alpha_2)<\kappa_{1\varepsilon}^0$, in the
parameter plane $(\alpha_1,\alpha_2)$ near the origin there exist
two curves:
${l}_h^\varepsilon$ and ${l}_\infty^\varepsilon$\\
\indent 1.\ \  the curve ${l}_h^\varepsilon$, which has the
following local representation:
$$
{l}_h^\varepsilon=\{(\alpha_1,\alpha_2):\
\kappa_2^\varepsilon(\alpha_1,\alpha_2)-\dfrac{b^\varepsilon}{a^\varepsilon}\kappa_1^\varepsilon(\alpha_1,\alpha_2)+\kappa_1^\varepsilon+h.o.t.=0,\
0<\kappa_1^\varepsilon(\alpha_1,\alpha_2)<\kappa_{1\varepsilon}^0\},
$$
is a Neimark-Sacker point branch of the numerical scheme (\ref{add}), where $h.o.t.=\varepsilon^3\cdot o(|(\alpha_1,\alpha_2)|)$, i.e.
$l_h^\varepsilon$ consists of
Neimark-Sacker bifurcation points of (\ref{add});\\
\indent 2.\ \  the curve ${l}_\infty^\varepsilon$, which has the
following local representation:
\begin{equation}\label{e2.4.6}
{l}_\infty^\varepsilon=\{(\alpha_1,\alpha_2):\
h(\alpha_1,\alpha_2)+h.o.t.=0,\ 0<\kappa_1(\alpha_1,\alpha_2)<\kappa_{1\varepsilon}^0\},
\end{equation}
is a homoclinic curve of the numerical scheme
(\ref{add}), where
$h(\alpha_1,\alpha_2)=\kappa_2^\varepsilon-\frac{6}{7}\frac{b^\varepsilon}{a^\varepsilon}\kappa_1^\varepsilon+\frac{5}{7}\kappa_1^\varepsilon$, $h.o.t.=\varepsilon^3\cdot o(|(\alpha_1,\alpha_2)|)$. In other
words, the numerical scheme (\ref{add}) presents  a unique homoclinic curve
connecting the origin for each
$(\alpha_1,\alpha_2)\linebreak[0]\in{l}_\infty^\varepsilon$. }
\end{theorem}

Comparing Theorem \ref{t44} to Theorem \ref{t2.4.1}, incorporating Theorem  \ref{t41} we obtain the following result.

\begin{theorem}\label{maintheorem}
Assume the assumption {\bf A} holds, that is the DDE (\ref{ep1}) exhibits a Takens-Bogdanov bifurcation at $(z,\alpha)=(0,0)$. Then the Takens-Bogdanov point of (\ref{ep1}) is inherited without any shift by the forward Euler scheme (\ref{add}) and turns into a 1:1 resonance point. Moreover,  there exists an $\varepsilon_0>0$, such that as $\varepsilon<\varepsilon_0$, the forward Euler scheme (\ref{add}) will reproduce the Hopf point branch and the homoclinic branch of the DDE (\ref{ep1}) with a shift of $O(\varepsilon)$ in parameter plane $(\alpha_1,\alpha_2)$, specially we have
$$|{l}_h-{l}_h^\varepsilon|=O(\varepsilon),$$
$$|{l}_\infty-{l}_\infty^\varepsilon|=O(\varepsilon).$$
\end{theorem}

\section{Numerical example}
In this section we present a numerical experiment to illustrate the theoretical results.

We consider a 1-dimensional DDE as follows
\begin{equation}\label{exam1}
\dot z(t)=(1+\alpha_1)z(t)-(1+\alpha_2)z(t-1)+\frac{1}{2}z(t)z(t-1).
\end{equation}
It is easy to show $(z,\alpha)=(0,0)$ is a Takens-Bogdanov point of (\ref{exam1}), cf. \cite{Fa95b,xu08}.
The forward Euler method for solving it is given by
\begin{equation}\label{numexam}
z_{k+1}=z_{k}+\varepsilon(1+\alpha_1)z_k-\varepsilon(1+\alpha_2)z_{k-m}+\frac{1}{2}\varepsilon z_kz_{k-m}
\end{equation}
with the step size $\varepsilon=\frac{1}{m}$, $m\in \mathbb{Z}_+$. Theorem \ref{t41} shows $(z,\alpha)=(0,0)$ is a 1:1 resonance point of (\ref{numexam}). The numerical experiment is carried for $m=100$, that is $\varepsilon=\frac{1}{100}$.

\begin{figure}
\centering
\includegraphics[width=10cm]{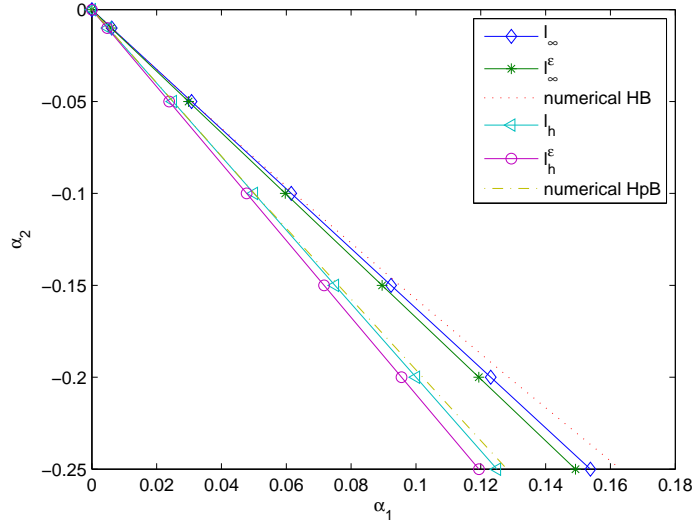}
\caption{Local bifurcation diagram in parameter plane $(\alpha_1,\alpha_2)$: by neglecting the h.o.t., the homoclinic branch $l_\infty$ and  Hopf point branch $l_h$ of (\ref{exam1}), the theoretical homoclinic branch $l_\infty^\varepsilon$ and Neimark-Sacker point branch $l_h^\varepsilon$ of the forward Euler discretization (\ref{numexam}), as well as the homoclinic branch "numerical HB" and neimark-Sacker point branch "numerical HpB" detected in Forward Euler discretization (\ref{numexam}) are plotted.}\label{fig1}
\end{figure}

From Theorem \ref{t2.4.1}, we know that $l_h=\{(\alpha_1,\alpha_2):\frac{4}{3}\alpha_1+\frac{2}{3}\alpha_2+h.o.t.=0\}$ is the local representation of the Hopf point branch, while $l_\infty=\{(\alpha_1,\alpha_2): \frac{26}{21}\alpha_1+\frac{16}{21}\alpha_2+h.o.t.=0\}$ is the local representation of the homoclinic branch of (\ref{exam1}) in parameter plane $(\alpha_1,\alpha_2)$. They are, by neglecting the higher order terms,  plotted in \textsc{Figure} \ref{fig1} by solid lines with triangles and diamonds, respectively.

From Theorem \ref{t44} we know that $l_h^\varepsilon=\{(\alpha_1,\alpha_2):(\frac{4}{3}+2\varepsilon)\alpha_1+(\frac{2}{3}-2\varepsilon)\alpha_2+h.o.t.=0\}$ is the local representation of the Neimark-Sacker point branch, while $l_\infty^\varepsilon=\{(\alpha_1,\alpha_2): (\frac{26}{21}+\frac{10}{7}\varepsilon)\alpha_1+(\frac{16}{21}-\frac{10}{7}\varepsilon)\alpha_2+h.o.t.=0\}$ is the local representation of the homoclinic branch of the forward Euler method (\ref{numexam}) in parameter plane $(\alpha_1,\alpha_2)$. They are, by neglecting the higher order terms, plotted in \textsc{Figure} \ref{fig1} by solid lines with circles and stars, respectively.

Besides, the realistic Hopf point branch of the forward Euler scheme (\ref{numexam}) is obtained by detecting the occurrence of the eigenvalues with modulus 1 at the fixed point $z=2(\alpha_2-\alpha_1)$. It is plotted in {\sc Figure} \ref{fig1} by dash-dot line. The critical values of parameter $\alpha$ for the occurrence of homoclinic orbits of the forward Euler scheme (\ref{numexam}) are obtained by a shooting technique. See the dotted line in {\sc Figure} \ref{fig1}.

\begin{figure}\centering
\subfloat[$\alpha=(0.005,-0.05)$]{\label{fig:subfig1}\includegraphics[width=4cm]{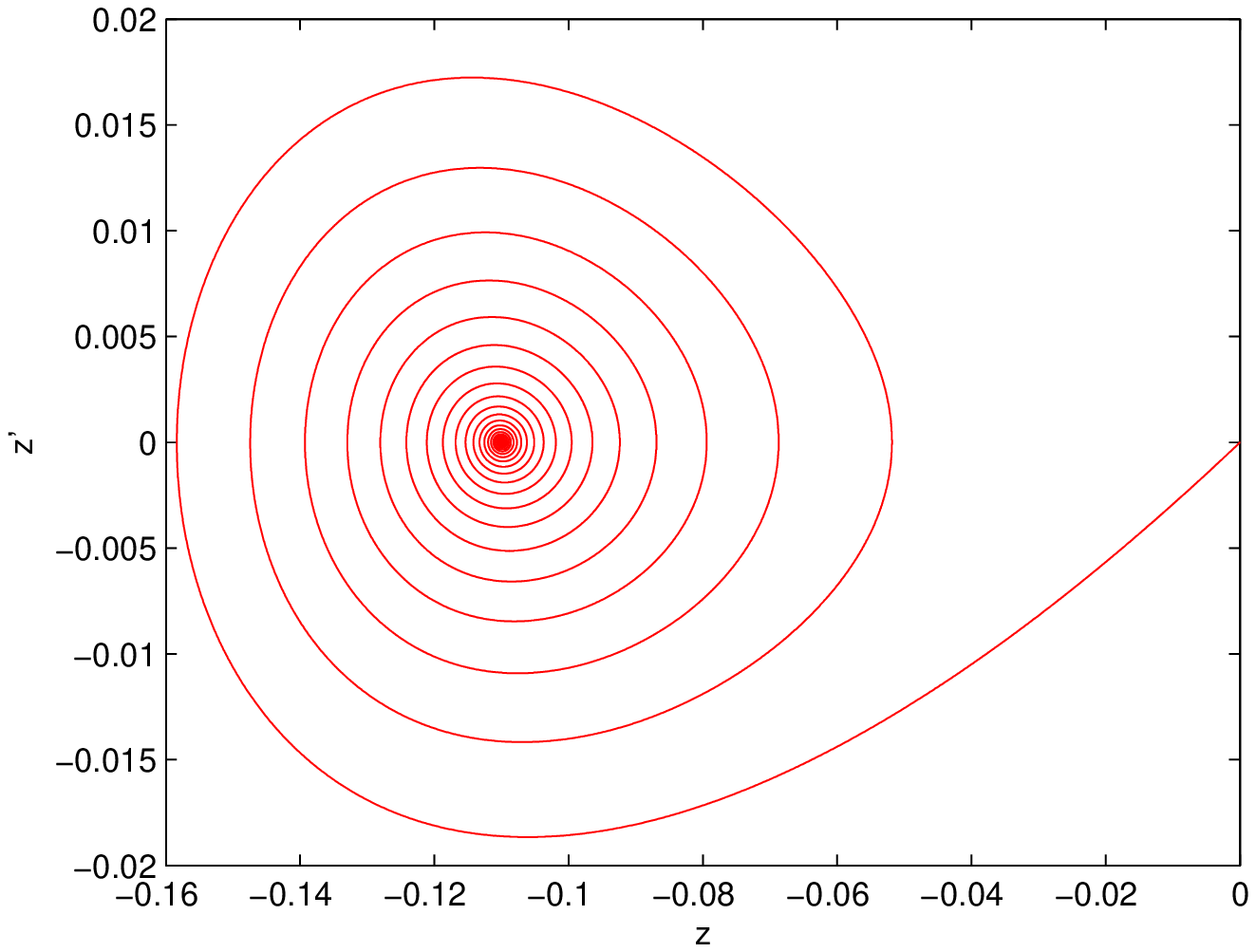}}
\subfloat[$\alpha=(0.05,-0.15)$]{\label{fig:subfig2}\includegraphics[width=4cm]{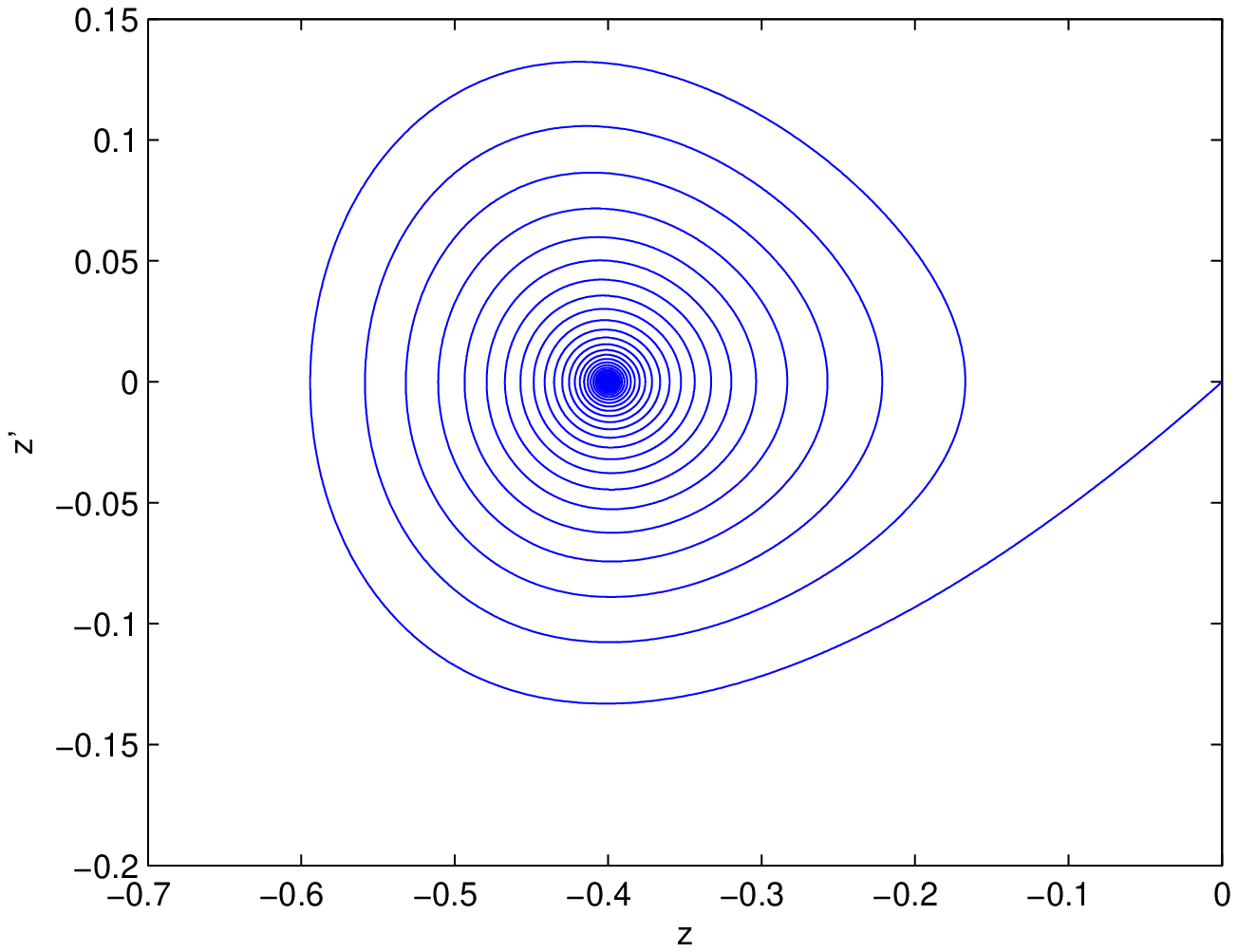}}
\subfloat[$\alpha=(0.08,-0.25)$]{\label{subfig3}\includegraphics[width=4cm]{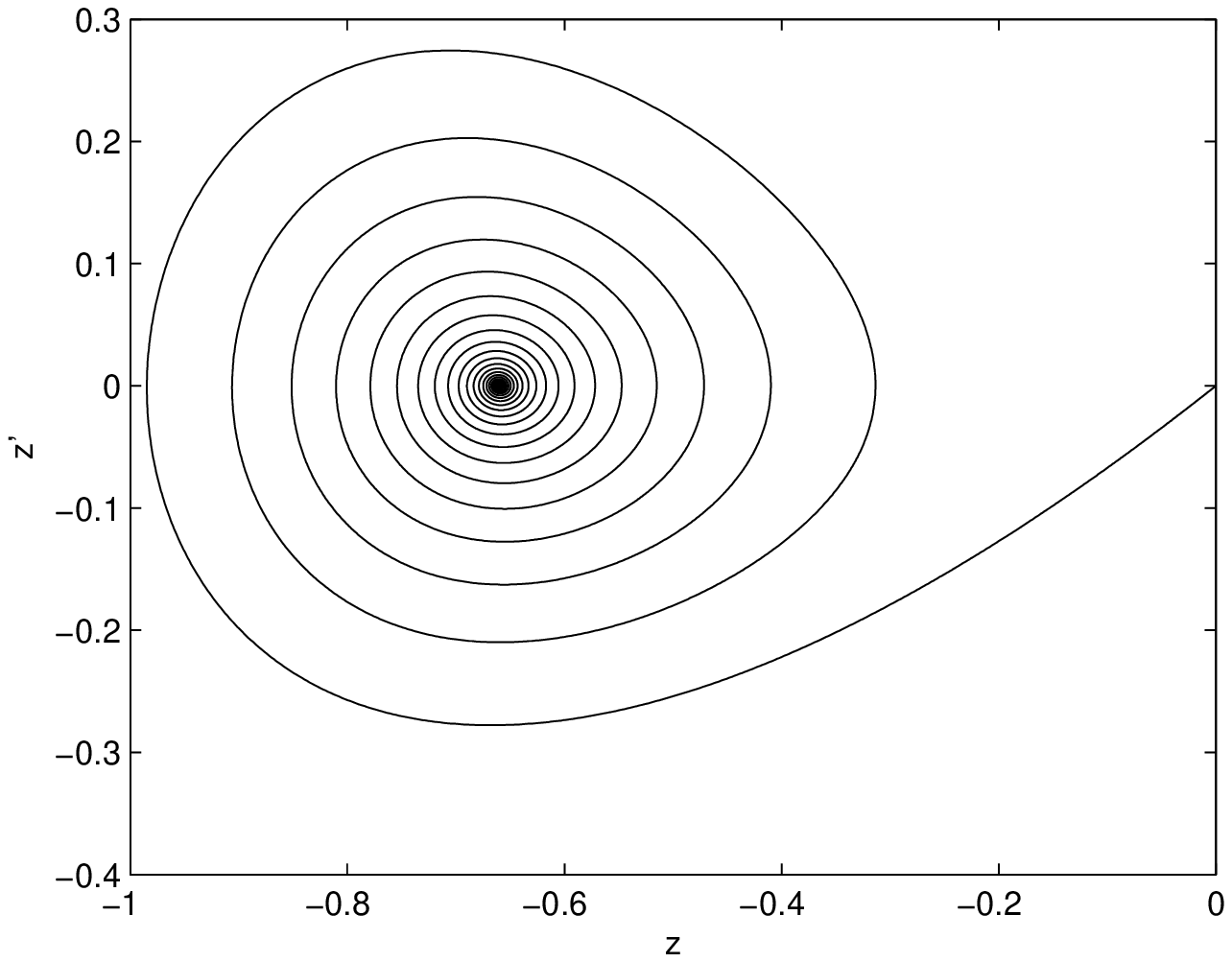}}
\caption{$2(\alpha_2-\alpha_1)$ are focuses of the Euler method (\ref{numexam}) when $\alpha$ belongs to the left region of $l_h^\varepsilon$. }\label{fig2}
\end{figure}

\begin{figure}\centering
\subfloat[$\alpha=(0.028,-0.05)$]{\includegraphics[width=4cm]{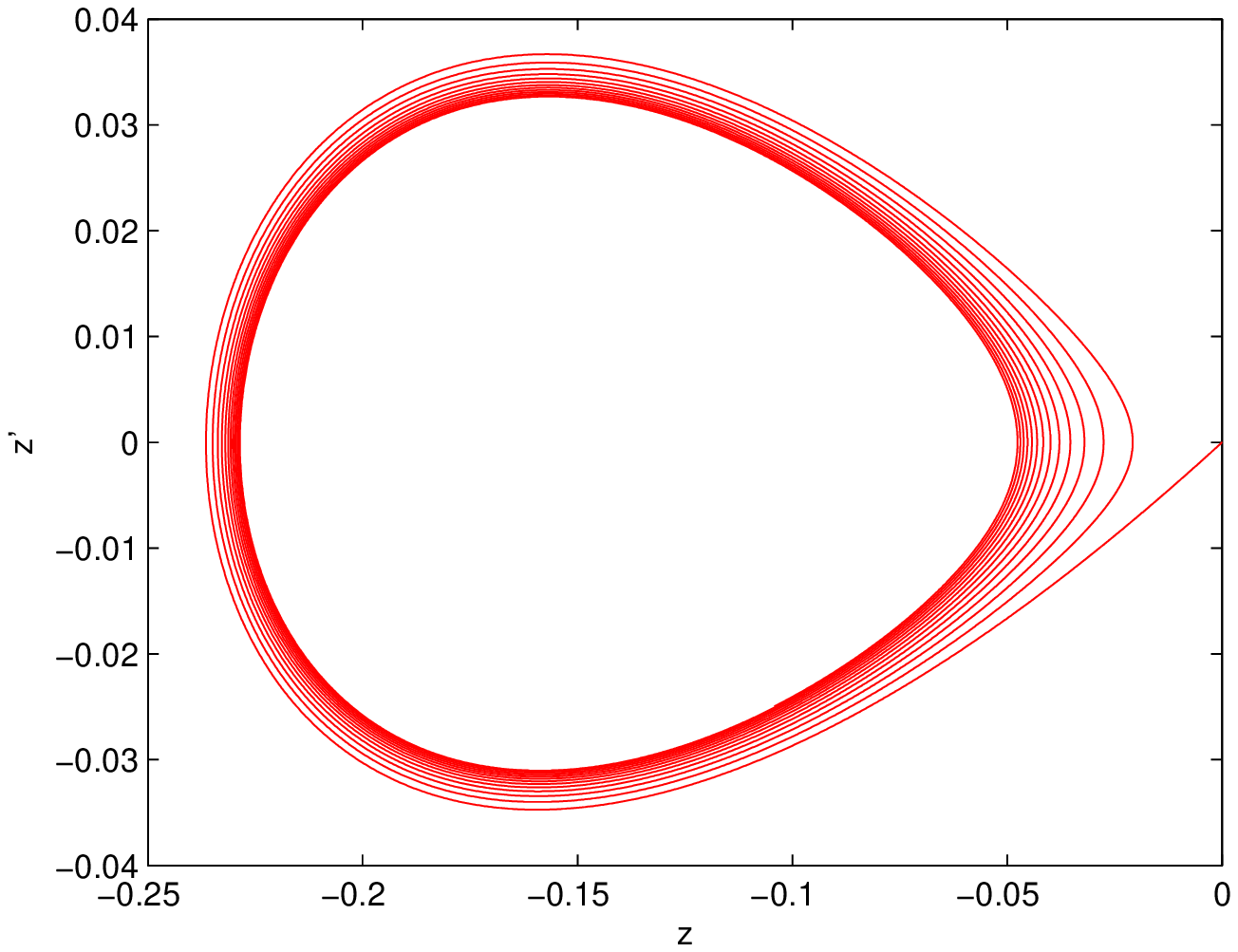}}
\subfloat[$\alpha=(0.085,-0.15)$]{\includegraphics[width=4cm]{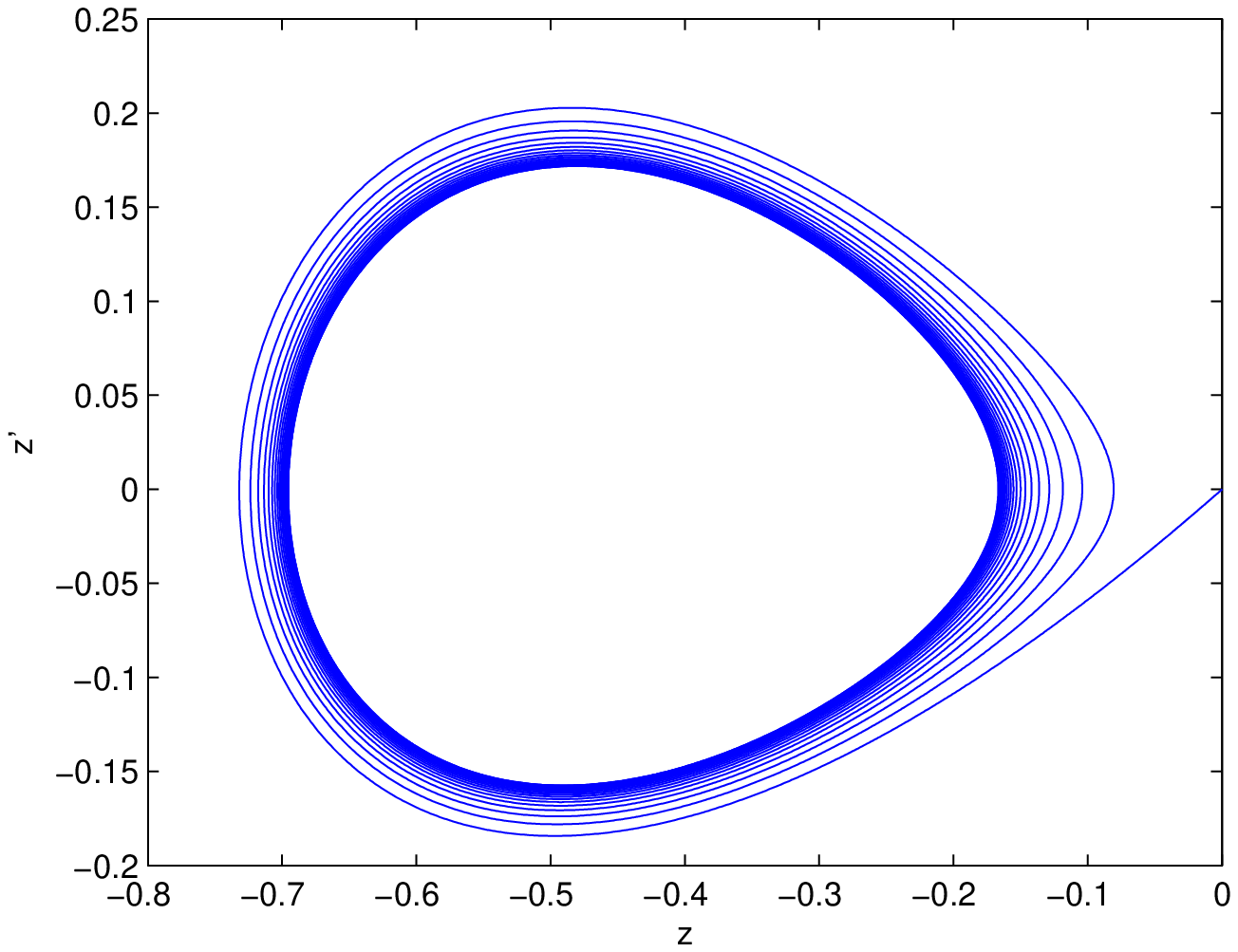}}
\subfloat[$\alpha=(0.145,-0.25)$]{\includegraphics[width=4cm]{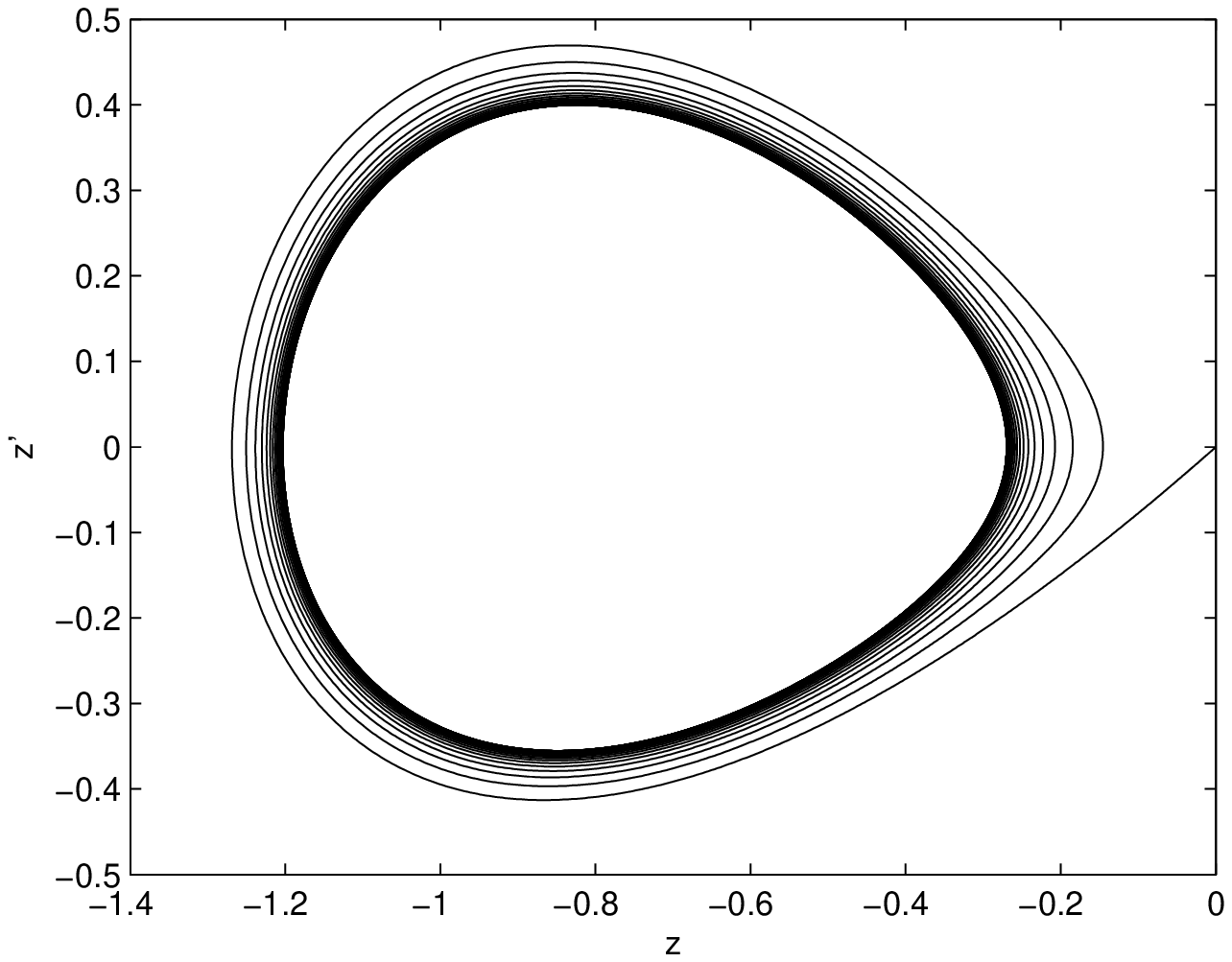}}
\caption{Periodic solutions bifurcated from $2(\alpha_2-\alpha_1$) of the Euler method (\ref{numexam}) when $\alpha$ is located in the region confined by $l_h^\varepsilon$ and $l_\infty^\varepsilon$.}\label{fig3}
\end{figure}

\begin{figure}\centering\includegraphics[width=10cm]{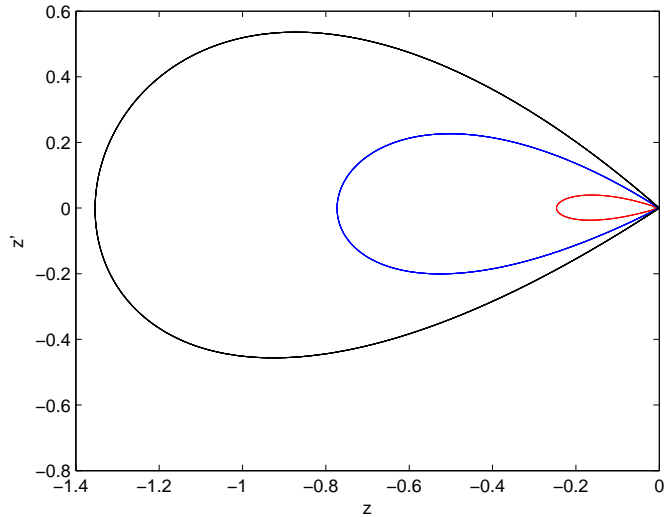}
\caption{Homoclinic orbits of (\ref{numexam}): $\alpha=(0.0308,-0.05)$ gives the small one, $\alpha=(0.0950,-0.15)$ leads to the middle one, and $\alpha=(0.1631,-0.25)$ results in the large one.}\label{fig4}
\end{figure}

In the parameter plane $(\alpha_1,\alpha_2)$, the fixed point of $z=2(\alpha_2-\alpha_1)$ is a focus when $(\alpha_1,\alpha_2)$ belongs to the left region of $l_h^\varepsilon$, when the parameter $(\alpha_1,\alpha_2)$ moves right and crosses $l_h^\varepsilon$, it turns into a central point and there will be periodic solutions bifurcating from this point. The periodicity will tend to infinity as the parameter $\alpha$ moves right and tends to $l_\infty^\varepsilon$. At last, the periodic solution becomes the homoclinic solution when $\alpha$ arrives at $l_\infty^\varepsilon$. These are the reproduction of the bifurcation structures for (\ref{exam1}) near the Takens-Bogdanov point. These processes are shown in {\sc Figure} \ref{fig2} to \ref{fig4} for fixed $\alpha_2=-0.05,-0.15,-0.25$, the corresponding values of $\alpha_1$ are $0.005, 0.05, 0.08$ (focuses), $0.028, 0.085, 0.145$ (periodic solutions) and $0.0308, 0.0950, 0.1631$ (homoclinic solutions), where the derivatives of $z$ are approximated by the difference quotient, and then a phase portrait like ODE's is presented.




\medskip

\end{document}